\documentstyle[12pt]{amsart}
\theoremstyle{plain}

\newtheorem{lem}{Lemma}

\newtheorem*{prop 1}{Proposition 1}
\newtheorem*{prop 2}{Proposition 2}
\newtheorem*{prop 3}{Proposition 3}
\newtheorem*{thm 1}{Theorem 1}
\newtheorem*{thm 2}{Theorem 2}
\newtheorem*{thm 3}{Theorem 3}
\newtheorem*{cond 1}{Condition a-s(1)}
\newtheorem*{cond 2}{Condition a-s(2)}
\begin{document}
\title[Presentations of symbolic dynamical systems]{Presentations of symbolic dynamical\\ systems by labelled directed graphs \\
(Notes for a "mini-cours", SDA2,\\ Paris 4-5 October 2007) }
\author{Wolfgang Krieger}
\begin{abstract}
We develop some aspects of a general theory of presentations of subshifts by labelled directed graphs, in particular  by compact graphs. Also considered are synchronization properties of subshifts that lead to presentations by countable graphs.
\end{abstract}
\maketitle
\def\det{{{\operatorname{det}}}}
\def\trace{{{\operatorname{trace}}}}
\def\card{{{\operatorname{card}}}}
Let $\Sigma$ be a finite alphabet.
In symbolic dynamics one studies subshifts $(X,S_{X})$, $X$ a 
shift-invariant closed subset of the shift space $\Sigma^{\Bbb Z}$  and $S_{X}$ 
the restriction of the left shift on  $\Sigma^{\Bbb Z}$    
to $X$. An introduction to 
the theory of subshifts is given in \cite {Ki} and in \cite {LM}. See here also 
\cite {BP}.

The first talk is about some aspects of a general theory of presentations of subshifts by labelled directed graphs, in particular  by compact graphs.  
The topic of the the second talk are synchronization properties of subshifts that lead to presentations by countable graphs.

We fix terminology and notation. Given a subshift $X \subset \Sigma^{\Bbb Z}$ we set
$$
	x_{[i,k]}Ê=Ê(x_{j})_{i \leq j \leq k}, \quad  		x \in X,  
			{i,k\in \Bbb Z} ,  i \leq k,
$$
and
$$
X_{[i,k]}Ê=Ê\{ x_{[i,k]} : x\in X \},\quad i,k\in{ \Bbb Z} ,  i \leq k . 
$$
We use similar notation also for blocks,
$$
b_{[i',k']}Ê=Ê(b_{j})_{i' \leq j \leq k'}, \quad 	b\in X_{[i,k]}, \; i \leq i' \leq k' \leq k,
$$
and also if indices range in semi-infinite intervals. ${\cal L} (X)$ is the set of admissible words of $X$. When convenient 
we identify blocks with the words they carry. We set
 \begin{align*}
&\Gamma_{n}^{+}(a)Ê= \{ b\in X_{(k,k+n]}: (a,b)\in X_{[i,k+n]}\} 
,\qquad n 
\in 
\Bbb N ,\\
&\Gamma_{\infty}^{+}(a)Ê= \{  y^{+}\in X_{(k,\infty )}:  (a, y^{+})\in 
X_{[i,\infty)}\} ,
\\
& \Gamma^{+}(a)Ê= \Gamma_{\infty}^{+}(a)Ê \cup \bigcup_{n\in \Bbb N}\Gamma^{+}_{n}(a)Ê, \qquad \qquad a\in X_{[i,k]},	 i,k\in{ \Bbb Z} ,  	i \leq k.
 \end{align*}
$ \Gamma^{-}$ has the time symmetric meaning. We set
 \begin{align*}
&\omega^{+}_{n}(a)Ê=
\bigcap_{x^{-} \in \Gamma_{\infty}^{-}(a)}
\{  b  \in X_{(k,k+n]}: (x^{-}, a,  b ) \in X_{( - \infty , k+n]}\},
\\
&\omega^{+}_{\infty}(a)Ê=Ê
\bigcap_{x^{-} \in \Gamma_{\infty}^{-}(a)}
\{ y^{+} \in X_{(k,\infty)}: (x^{-}, a, y^{+}) \in X\},
\\
& \omega^{+}(a)Ê= \omega_{\infty}^{+}(a)Ê \cup \bigcup_{n\in \Bbb N}\omega^{+}_{n}(a)Ê, \qquad a\in X_{[i,k]},	 i,k\in{ \Bbb Z} ,  	i \leq k.
 \end{align*}	
$\omega^{-}$ has the time symmetric meaning. 
The unstable set of a point $x \in X$ is denoted by $W^{-}_{X}(x)$,
$$
W^{-}_{X}(x) = \bigcup_{I \in{ \Bbb N}} \{ y\in X: y_{i}= x_{i}, i \leq I\}.
$$
The $W^{-}_{X}(x), x \in X,$ carry the inductive limit topologies of the compact topologies on the sets 
$\{ y\in X: y_{i}= x_{i}, i \leq I\}, I \in \Bbb N  $.

We recall that, given subshifts $X \subset \Sigma^{\Bbb Z},  \bar{X}\subset \bar{\Sigma}^{\Bbb Z}$,
and a  topological conjugacy $\varphi: X \rightarrow \bar{X} 
$, 
there is 
for some 
$L\in{ \Bbb Z}_{+} $ a block mapping
$$
\Phi:  X_{[-L,L]} \rightarrow \bar{\Sigma } 
$$
such that
$$
\varphi (x)Ê=Ê(\Phi (x_{[i-L,i+L]}))_{i\in {\Bbb Z}} .
$$
We say then that $\varphi$ is given by $\Phi$, and we write
$$
\Phi( a)Ê=Ê(\Phi (a_{[j-L,j+L]}))_{i+L\leq j\leq k-L},	\quad a\in 
X_{[i,k]},  \quad i,k\in{ \Bbb Z} ,	k - i \geq 2L,
$$
and use similar notation if indices range in semi-infinite intervals. The interval $[-L,L]$ is called a coding window.
\section{Presentations of subshifts by Shannon graphs}
Let $\Sigma$ be a finite alphabet and consider a directed graph ${\mathcal G}$ with
vertex set ${\mathcal V}$ whose directed edges are labeled with symbols in 
$\Sigma$.  The graph ${\mathcal G}$ is called Shannon if its labeling is 1-right
resolving in the sense that for every vertex $V \in{\mathcal V} $ and for every
symbol $\sigma \in \Sigma $ there is at most one edge in ${\mathcal G}$ that leaves
 $V$ and that carries the label $\sigma$. Shannon graphs are 
 also known as
deterministic transition systems. Denote the set
of vertices $V \in{ \mathcal V}$ of a Shannon graph ${\mathcal G}$ that have an outgoing
edge that carries the label $\sigma \in \Sigma$ by $\cal V(\sigma)$, and for
$V \in{ \mathcal V}(\sigma) $ denote by $\tau_{\sigma}(V)$ the final vertex of the
edge that
leaves $V$ and that carries the label $\sigma$. We call $(\tau_{\sigma})
_{\sigma \in \Sigma}$ the transition rule of the Shannon graph. The forward
context $\Gamma^{+}_{\infty}(V)$ of a vertex $V$ of a Shannon
graph ${\mathcal G}$ is defined as the closed subset of $\Sigma^{\Bbb N}$ that
contains the label sequences of the semi-infinite paths in ${\mathcal G}$ that leave
$V$. A Shannon graph is called forward separated if distinct vertices have
distinct forward contexts. Call a Shannon graph compact if its vertex set
$\mathcal V $ carries a compact topology such that the sets ${\mathcal} V(\sigma)
$, are open and such that the mappings 
$V\to\tau_{\sigma}(V), (V\in{\mathcal V}(\sigma)),\sigma \in \Sigma$, are continuous.
We say that a Shannon graph ${\mathcal G}$, in which
every vertex has at least one
incoming edge and at least one outgoing edge presents a subshift $X \subset
\Sigma^{\Bbb Z}$, if the set of admissible words of $X$ coincides with the set
of label sequences of the finite paths in $\mathcal G$.

For a finite alphabet $\Sigma$ denote by ${\mathcal V}(\Sigma)$ the set of non-empty closed
subsets of $\Sigma^{\Bbb N}$ with its compact Hausdorff subset topology. ${\mathcal
V}(\Sigma)$ is the vertex set of a compact Shannon graph $\mathcal G( \Sigma)$: For
$\sigma \in \Sigma$ the set ${ \mathcal V} (\Sigma) (\sigma)$ of ${\mathcal G}(\Sigma)$
is equal to the set of $V \in{\mathcal V}$ that contain a sequence that 
starts with $\sigma$,
and the transition rule for ${\mathcal G}(\Sigma)$ is
 \begin{align*}
\tau_{\sigma}(V) = \{v_{(1, \infty)}: v \in  V, v_{1} = \sigma \}, \qquad
V\in
{\mathcal V}(\Sigma)(\sigma ), \sigma \in \Sigma. 
\tag{1.1}
 \end{align*}
The sets  ${\mathcal V}(\Sigma)(\sigma), \sigma \in \Sigma$, are compact-open and we can
associate to
the finite alphabet $\Sigma$ also the topological Markov chain 
$tM(\Sigma)$
which is the compact dynamical system that is obtained by having the left shift
act on the space
$$
\{ (x_{i},V_{i})_{i\in {\Bbb Z}} \in (\Sigma \times {\mathcal V} (\Sigma))^{\Bbb Z} :
V _{i+1}   =
\tau_{x_{i+1}}(V _{i} ), i\in{ \Bbb Z}\}.
$$

Denote by ${\cal C} (W^-_{\Sigma^{\Bbb Z}}(x))$ the set of closed subsets of
$W_{\Sigma^{\Bbb Z}}^{-}(x), x \in \Sigma^{\Bbb Z}$.
For the finite alphabet $\Sigma$ we introduce the set 
$ {\mathcal C}^{\bullet}(\Sigma)$ of pairs $(x,C)$, where $x \in \Sigma^{\Bbb Z}$ and where
$C\in{\mathcal C}
 (W^{-}(x)) $. 
Having the left shift act on $ { \mathcal 
C}^{\bullet}(\Sigma)$
one obtains a dynamical system for which one has a shift commuting bijection
$\psi_{\Sigma}$ onto $tM(\Sigma)$. This bijection $\psi_{\Sigma}$ assigns to a
point $(x,C) \in { \mathcal C}^{\bullet}(\Sigma)$ the point   $(x_{i}, V_{i}  
(C))_{i\in
{\Bbb Z}} \in tM(\Sigma)$ that is given by
$$
 V_{i}(C)   = \{ y \in C: x_{j} = y_{j} , j \leq i \}, \qquad i \in {\Bbb Z}.
$$
In order to turn  ${ \mathcal C}^{\bullet}(\Sigma)$  into a compact dynamical system one transports the
topology on
$tM (\Sigma)$ to $ {\mathcal C}^{\bullet}(\Sigma)$   by means  of  the inverse of the bijection
$\psi_{\Sigma}$.
\begin{lem}
Let $X \subset \Sigma^{\Bbb Z}, \widetilde{X }\subset
\widetilde{\Sigma}^{\Bbb Z}$ be subshifts, and let $\varphi: X \to \widetilde{X
}$
be a topological conjugacy. 
Let $L \in { \Bbb Z}_{+} $ be such that $[-L,L]$ is a coding window for 
$\varphi$ and for $\varphi^{-1}$. Let $\varphi$ be given by the block map
$\Phi: X_{[-L, L]}\to \widetilde{\Sigma},$ and  let $\varphi^{-1}$  be given  by the block map
$ \widetilde{\Phi}:  \widetilde{X}_{[-L, L]}\to \Sigma.$ 
Let $(x, C) \in {\mathcal C}^\bullet(\Sigma)$, and set $  ( \widetilde{x},  \widetilde{C})  = \varphi(x, C)$.
Then
\begin{multline*}
V_i(  \widetilde{C}) = \{  \widetilde{\Phi}( y^+ ):   y^+\in V_{i-L}(C), \Phi( x_{[i-3L, i-L]},y^+_{[i-L, L]} )
=    \widetilde{x}_{[i-2L, i]}  \},
\\ i \in{ \Bbb Z}.
\end{multline*}
\end{lem}
\begin{pf}
One notes that
\begin{multline*}
V_i(  \widetilde{C})  \supset \{  \Phi(x^+_{[ i-3L, i+L   ]} y^+  ) : y^+ \in V_{i+L}(C)     \}, \quad
 \widetilde{\Phi}(V_i( \widetilde{C}   ) \subset V_{i-L}( C) ,\\ i \in \Bbb Z . \qed
\end{multline*}
\renewcommand{\qedsymbol}{}
\end{pf}
We call a subset ${\mathcal V}$ of ${\mathcal V} (\Sigma)$ transition-complete if for $V
\in  {\mathcal V} \cap {\mathcal V} (\Sigma)(\sigma)$ also $\tau_{\sigma}(V) \in {\mathcal V}$. To a transition complete subset $\cal V$ of $\cal V (\Sigma)$ we associate
 the sub-Shannon graph ${\mathcal} G(\cal V)$ of  ${\mathcal G}(\Sigma)$  that has as
vertex set the set 
${\mathcal V}$ and as transition rule the restriction of the rule (1.1) to  ${\mathcal 
V}$. The set of forward contexts of the vertices af a forward separated Shannon 
graph is a transition-complete subset of ${\mathcal V} (\Sigma)$, and
the mapping that sends every vertex of a forward separated Shannon graph 
to its
forward context is an isomorphism of ${\mathcal G}$ onto the Shannon graph that 
is associated with the set of forward contexts of ${\mathcal G}$.

We say that a transition-complete subset ${\mathcal V}$ of  $  \cal V 
(\Sigma)$ is
retro-complete if every $V \in \cal V$ thas a a predecessor in ${\mathcal V} $. For
a transition-complete compact subset ${\mathcal V} $ of $ \cal
V (\Sigma)$ the
sets $ {\mathcal V}(\sigma) , \sigma \in \Sigma$, are also compact, and we set
$$
\tau ( {\mathcal V} ) = \bigcup _{\sigma \in \Sigma} \tau_{\sigma} ( {\mathcal V} (\sigma)),
$$
and with $\tau^{(0)} ({\mathcal V}) = {\mathcal V} $, we set inductively
$$
\tau^{(n)} ({\mathcal V} ) =  \tau(\tau^{(n-1)}({\mathcal V})), \qquad n \in \Bbb N. 
$$
Here
$$
\tau^{(n)} {\mathcal V}) \subset   \tau^{(n-1)}({\mathcal V}),  \qquad n \in {\Bbb N}, 
$$
and the intersection $\bigcap_{n \in \Bbb N} \tau^{(n)} ({\mathcal V})  $ is the
maximal transition-  and retro-complete subset of ${\mathcal V}$.

The transition- and retro-complete subsets of $\cal V(\Sigma)$ are in
one-to-one correspondence with the shift invariant subsets of $tM(\Sigma)$: To
a transition- and retro-complete set ${\mathcal V} \subset  \cal V(\Sigma)$ there
corresponds the system
$$
tM({\mathcal V}) = \{ (x_{i},V_{i} )_{i\in {\Bbb Z}} \in (\Sigma \times {\mathcal V} )^{{\Bbb
Z}} : V
_{i+1}   =
\tau_{x_{i+1}}(V _{i} ), i\in {\Bbb Z}\}.
$$ 
Also, assigning to a transition- and retro-complete subset $\cal V$ of 
$\{{\mathcal V}(\Sigma)$  the dynamical system $C^ {\bullet}({\mathcal V} )=
\psi^{-1}_{\Sigma}(tM({\mathcal V}))$, sets up a one-to one corresondence between the
transition- and retro-complete subsets of ${\mathcal V}(\Sigma)$ and the shift
invariant subsystems $C^ {\bullet}$ of
$C^ {\bullet}(\Sigma)$ with the property that $(x, C) \in  C^ {\bullet}$  and
$y \in   C$ imply hat $(y,C) \in  C^ {\bullet} $.

\begin{lem}
Let $X \subset \Sigma^{\Bbb Z}$ be a subshift. Let $x^{-} \in X_{(-\infty,
0]}, C \in {\mathcal V}(\Sigma$), and 
$$
x^{-}(k) \in X_{(-\infty, 0]},\quad  C(k) \in{\mathcal V}(\Sigma), C(k) \subset
\Gamma^{+}_{\infty}(x^{-}(k) ), \qquad k \in \Bbb N,
$$
and let
$$
x^{-} = \lim _{k \to \infty} x^{-}(k) ,\quad C = \lim _{k \to \infty} C(k). 
$$
Then 
 \begin{align*}
C \subset  \Gamma^{+}_{\infty}(x^{-})  .   \tag {1.2}
 \end{align*}
\end{lem}
\begin{pf}
For $n\in{ \Bbb N}$ let  $k_{n}  \in{ \Bbb N} $ be such that
$$
x^{-}_{[-n,0]} =x^{-}_{[-n,0]} (k), \quad C_{[1,n]} = 
C_{[1,n]}(k),\qquad k \geq k_{n}.
$$
Then for $x^{+} \in C$ and $ n\in {\Bbb N} $
$$
(x^{-}_{[-n,0]},x^{+}_{[1,n]}) = (x^{-}_{[-n,0]}(k),x^{+}_{[1,n]}) \in X_{[-n,
n]},\qquad k \geq k_{n} ,
$$
which implies (1.2).\hfill \qed
\renewcommand{\qedsymbol}{}
\end{pf}
\begin{prop 1}
Let $X \subset \Sigma^{\Bbb Z}$ be a subshift. The set
$$
{\mathcal V}_{\circ}(X) = \bigcup_{x^{-} \in X_{(- \infty, 0]} }\{ V \in {\mathcal V}
(\Sigma): V \subset \Gamma^{+}_{\infty}( x^{-}) \}
$$
is closed.
\end{prop 1}
\begin{pf}
Apply Lemma 2 and the compactness of $X_{(- \infty, 0]}$. \qed
\renewcommand{\qedsymbol}{}
\end{pf}
For a subshift $X \subset \Sigma^{{\Bbb Z}}$ we set
$$
{\mathcal V} _{max}(X) = \bigcap_{n\in {\Bbb N}} \tau^{n}({\mathcal V}_{\circ}(X)),
$$
${\mathcal G}({\mathcal V} _{max}(X))$ presents $X$ and we refer to ${\mathcal G}({\mathcal V}
_{max}(X))$ as the maximal presenting Shannon graph of
$X$. This terminology is justified by the fact that every transition-
and retro-complete subset ${\mathcal V}$ of ${\mathcal V} (\Sigma)$, such that ${\mathcal G}({\mathcal V})  
 $ presents $X$, is a subset of ${\mathcal V} _{max}(X)$.
\begin{prop 2}
Let $X \subset \Sigma^{{\Bbb Z}}, \widetilde{X }\subset
\widetilde{\Sigma}^{{\Bbb Z}}$ be subshifts, and let $\varphi: X \to \widetilde{X}$
be a topological conjugacy. Then
$$
(x,C) \to \varphi(x,C) \quad ((x,C)   \in C^{\bullet}({ \mathcal V}
_{max}(X)))
$$
is a topological conjugacy of $C^{\bullet}( { \mathcal V} _{max}(X))$ onto
$C^{\bullet}({ \mathcal V} _{max}(\widetilde{X })$.
\end{prop 2}
\begin{pf}
An application of Lemma 1 yields the continuity of the mapping
$$
(x,C) \to \varphi(x,C) \quad ((x,C)   \in C^{\bullet}( { \mathcal V}_{max}(X)). \qed
$$
\renewcommand{\qedsymbol}{}
\end{pf}

Given subshifts $ X \subset \Sigma^{{\Bbb Z}}, \widetilde{X }\subset
\widetilde{\Sigma}^{{\Bbb Z}}$ and a topological conjugacy $\varphi: X
\to \widetilde{X }$, $\varphi$ will also denote the topological conjugacy that
sends the point $(x,C)   \in C^{\bullet}({ \mathcal V} _{max}(X))$  to the point
$\varphi(x,C) \in C^{\bullet}({ \mathcal V} _{max}(\widetilde{X })$. 
Given
transition- and retro-complete subsets ${ \mathcal V} \subset{ \mathcal V}_{max}(X)$ and
$\widetilde{{ \mathcal V}}\subset{ \mathcal V}_{max}(\widetilde{X })$, where ${ \mathcal G} ({ \mathcal V}
)$ presents $X$ and ${ \mathcal G}(\widetilde{{ \mathcal V} })$ presents $\widetilde{X }$,
such that  $\varphi(C^{\bullet}({ \mathcal V} )) = C^{\bullet}(\widetilde{{ \mathcal V} })$, we
write also $\varphi({ \mathcal V}) = \widetilde{{ \mathcal V} }$ and $\varphi(tM({ \mathcal V})) =
tM(\widetilde{{ \mathcal V}})$.

Given a construction for a subshift $X \subset \Sigma^{{\Bbb Z}}$  of a
transition- and retro-complete subset ${ \mathcal V} _{X} \subset{ \mathcal V}_{max}(X)$
such that ${\mathcal G}({ \mathcal V} _{X} )$ presents $X$, we say that ${ \mathcal V} _{X}$ is 
canonical, if for subshifts $
X \subset \Sigma^{{\Bbb Z}}, \widetilde{X }\subset
\widetilde{\Sigma}^{{\Bbb Z}}$  and a topological conjugacy $\varphi: X
\to \widetilde{X }$, $\varphi ({ \mathcal V} _{X}) = { \mathcal V} _{\widetilde{X }}$. 
${ \mathcal V} _{max}$ itself is canonical by Proposition 2. The
standard example of the canonical situation is  
$$
{ \mathcal V}_{standard}(X) = \{\Gamma^{+}_{\infty}(x^{-}): x^{-} \in X_{(- \infty,
0]} \}.
$$
It is
$$
tM({ \mathcal V}_{standard}(X) ) = \{(x_{i}, \Gamma^{+}_{\infty}((x_{j})_{j<i})
)_{i\in {\Bbb Z}}: (x_{i})_{i\in{ \Bbb Z}} \in X\},
$$
and
$$
C^{\bullet}({ \mathcal V}_{standard}(X)) = \{(x, W^{-}(X)): x \in X \},
$$
and from this it is seen that $ { \mathcal V}_{standard}$ is canonical. The finiteness of $  { \mathcal V}_{standard}$ characterizes the sofic case  \cite  {W}. That  $  { \mathcal V}_{standard}$ is canonical was first noted in the sofic case in  
\cite {Kr1, Kr2}.
For a subshift $X \subset \Sigma^{{\Bbb Z}}$ the closure of ${ \mathcal V}_{standard}(X)$  is a canonical compact 
Shannon graph that presents $X$. This presentation appears in Matsumoto's theory of 
$\lambda$
-graph systems (see \cite {M, KM}). Also in non-sofic cases
${ \mathcal V}_{standard}(X)$ itself can be compact. For instance, the coded system (see \cite {BH})
with alphabet $\Sigma = \{ \gamma, 0, 1\}$ and code $\{\gamma 0^{l}  0^{l} : l
\in \Bbb N \}$ has a compact  ${ \mathcal V}_{standard} $, while  the coded system
with  alphabet $\Sigma = \{ \gamma, 0, 1\}$ and code $\{\gamma 0^{2l} 
0^{2l} : l \in {\Bbb N} \}$ has a ${ \mathcal V}_{standard}$, that is not compact. The content of
the following
proposition is that the compactness of ${ \mathcal V}_{standard} $ is an invariant of
topological conjugacy.

\begin{prop 3}
Let $X \subset \Sigma^{\Bbb Z}$ and $\widetilde{X }\subset
\widetilde{\Sigma}^{\Bbb Z}$ be topologically conjugate subshifts,
and let ${ \mathcal V}_{standard}(X) $ be compact. Then ${ \mathcal V}_{standard}(\widetilde{X }) $ is also compact.
\end{prop 3}
\begin{pf}
This follows from Proposition 2.
\end{pf}

Another example of a canonical presenting Shannon graph  is the word Shannon graph of a subshift $X \subset 
\Sigma^{{\Bbb Z}}$,
$$
{\mathcal V} _{word}(X) = \{\{x^{+}\}: x^{+} \in X_{[1, \infty)}\}.
$$
It is
$$
tM({\mathcal V} _{word}(X) ) = \{(x_{i},\{ (x_{j})_{j \geq i}\})_{i \in{ \Bbb Z}}:
(x_{i})_{i \in{ \Bbb Z}} \in X \},
$$
and
$$
C^{\bullet}({\mathcal V}_{word}(X)) = \{(x, \{x\}): x \in X \}.
$$
${\mathcal V} _{word}$ is compact.

\section{Notions of Synchronization}
We describe synchronizing shifts and their synchronizing 
Shannon graphs and then consider  the more general notions of s-synchronization and 
a-synchronization.
\subsection{Synchronization}
A word $b$ that is admissible for a topologically transitive subshift  $X\subset \Sigma^{{\Bbb Z}}$ 
is called a synchronizing word of $X$, 
if for $c \in \Gamma^{-}(b), d
\in \Gamma^{+}(b) $ one has that  $cbd \in {\mathcal L} (X)$. Equivalently, a
synchronizing word of $X$ can be defined as a word $b\in  {\mathcal L}(X)$ such that
$\Gamma^{-}(b) =
\omega^{-}(b)$, or, such that  $\Gamma^{+}(b) = \omega^{+}(b)$. We denote the
set of synchronizing words  of a topologically transitive subshift $X \subset
\Sigma^{{\Bbb Z}}$ by ${\mathcal L}
_{synchro}(X)$. 
\begin{lem}
Let $X \subset \Sigma^{{\Bbb Z}}, \widetilde{X} \subset \widetilde{\Sigma}^{{\Bbb Z}}$ be topologically transitive subshifts and let $L\in {\Bbb Z}_{+}$ be such
that there is a topological conjugacy of $X$ onto $\widetilde{X}$ that has
together with its inverse the coding window $[-L, L]$, with its inverse  given by
a block map
$$
\widetilde {\Phi}: \widetilde{X}_{[-L, L]} \to \Sigma.
$$
Let  
$$
b \in {\mathcal L}_{synchro}(X),
$$
and let $\widetilde{b}\in {\mathcal L} (
\widetilde{X} ) $ be a word such that $b =\widetilde {\Phi} (\widetilde{b})$. Then
$$
\widetilde{b} \in{\mathcal L}_{synchro}(\widetilde{X}),
$$
\end{lem}
\begin{pf}
Let $\widetilde{c} \in \Gamma^{-}( \widetilde{b}),\widetilde{d} \in
\Gamma^{+}( \widetilde{b})$ and choose words $\widetilde{c}^{\prime} \in
\Gamma^{-}( \widetilde{c} \widetilde{b})$ and $ \widetilde{d}^{\prime}  \in
\Gamma^{+}(  \widetilde{b}\widetilde{d}) $ of length $2L$. With words $ c\in
\Gamma^{-}(b),d\in \Gamma^{+}(b),  $ that are given by
$$
\widetilde{\Phi}( \widetilde{c}^{\prime}\widetilde{c}\widetilde{b} ) = cb,
\quad  \widetilde{\Phi}( \widetilde{b}\widetilde{d}\widetilde{c}^{\prime}\ ) =bd,
$$
one has 
$$
cbd \in{\mathcal L} (X),
$$
and
$$
\widetilde{c} \widetilde{b} \widetilde{d}    =  \Phi  (cbd).
\qed
$$
\renewcommand{\qedsymbol}{}
\end{pf}
A topologically transitive subshift $X \subset \Sigma^{{\Bbb Z}}$ that has a synchronizing word is called synchronizing
 (see \cite {BH}).
For $a \in {\mathcal L} _{synchro}(X)$ and $b \in \Gamma ^{+}(a)$ also $ab \in {\mathcal L} _{synchro}(X)$. It follows for a synchronizing subshift $X$  that the set
$$
{\mathcal V} _{synchro}(X) = \{ \Gamma^{+}_{\infty}(a) : a \in{\mathcal L} _{synchro}(X) \}
$$
is the vertex set of an irreducible Shannon sub-graph ${\mathcal G}({\mathcal V}_{synchro}(X))$  of ${\mathcal G}({\mathcal V}_{max}(X))$, that we call
the synchronizing Shannon graph of $X$. ${\mathcal G}({\mathcal V}_{synchro}(X))$   presents $X$.
Topologically transitive sofic systems can be characterized as  the synchronizing subshifts whose synchronizing Shannon graph is finite 
\cite {W}. The Shannon graph
${\mathcal G}({\mathcal V}_{synchro}(X))$ was first constructed in the sofic case in 
\cite {F},
where the term "Shannon graph"  was introduced.
It follows from Lemma 3  that
synchronization is an invariant of topological conjugacy.
The presentation of a synchronizing subshift by its synchronizing Shannon graph
is canonical. This is the content of the following theorem. 
\begin{thm 1}
\label{thm 1}
Let $X \subset \Sigma^{\Bbb Z}, \widetilde{X} \subset \widetilde{\Sigma}^{\Bbb
Z}$  be synchronizing subshifts, and let $\varphi: X \to \widetilde{X}  $ be a
topological conjugacy. Then
$$
\varphi ({\mathcal V} _{synchro}(X)) =    {\mathcal V}_{synchro}(\widetilde{X} ).
$$
\end{thm 1}
\begin{pf}
Let $L\in{ \Bbb Z}_+$ be such that $[-L, L]$ is a coding window for $\varphi$  and for $\varphi ^{-1}$.

Let
$$
(x, C) \in   C^{\bullet }({\mathcal V} _{synchro}(X)).
$$
and let $j \in {\Bbb Z} $. Let $b$ be a synchronizing word of $X$ such that
\begin{align*}
V_{-3L}(C) = \Gamma^{+}_{\infty}(b),   \tag {2.1}
\end{align*}
end let $d \in \Gamma^{+}(b) $ be such that $bd \in \Gamma^{-}(b) $. Consider the point
 $x^{\prime} \in 
X $  such that  $X_{(- \infty, - 3L)} $ carries the left infinite concatenation
of $db$, and such that
\begin{align*}
x^{\prime}_{[-3L, \infty)}= x_{[-3L, \infty)}. \tag {2.2}
\end{align*}
Let $( x^{\prime} , C^{\prime}    ) \in C^{\bullet }({\mathcal V}_{synchro}(X))  $ be given by
$$
V_i( C^{\prime}   ) = \Gamma^+_\infty( x^{\prime}_{(-\infty  , i]}    ), \quad i\in{ \Bbb Z},
$$
and let
$$
(  \widetilde{x}^{\prime} , \widetilde{C}^{\prime}    ) = \varphi ( x^{\prime} , C^{\prime}    ).
$$
By Lemma 3 
$$
(  \widetilde{x}^{\prime} , \widetilde{C}^{\prime}    ) \in C^{\bullet }({\mathcal V}_{synchro}(\widetilde{X})),
$$
and by Lemma 1 and by (2.1) and (2.2)
$$
V_j( \widetilde{C}^{\prime}  ) =V_j( \widetilde{C}  ) .
$$
By symmetry the theorem is proved. \qed
\renewcommand{\qedsymbol}{}
\end{pf}
\subsection{ s-synchronization}
A word $b$ that is admissible for a topologically transitive subshift  $X\subset \Sigma^{{\Bbb Z}}$ 
is called an s-synchronizing word of $X$, if for all
$c \in {\mathcal L}(X)$ there exists a $ d \in \Gamma^{+}(c)$ such that $cd \in
\omega^{-} (b). $ A synchronizing word is s-synchronizing. We denote the set
of s-synchronizing words of a topologically transitive subshift $X \subset
\Sigma^{\Bbb Z}$ by ${\mathcal L}_{s-synchro}(X)$. We note that for an
s-synchronizing word $b$ there exists in particular a word $d \in
\Gamma^{+}(b) $ such that $bd \in \omega^{-} (b) $ and this implies that
$\omega^{-}_{\infty}(b) \neq \emptyset$.
\begin{lem}
Let $X 
\subset \Sigma^{\Bbb Z}$ be a topologically transitive subshift, and let 
\begin{align*}
b\in {\mathcal L}_{s-synchro} (X). \tag{2.3}
\end{align*}
Let $d \in{\mathcal L} (X)$ be
such that $bd \in \omega^{-}(b)$ and let $l$ denote the length
of the word $bd$. 

Let $\widetilde{X} \subset \widetilde{\Sigma}   ^{\Bbb Z}$ be a
subshift that is topologically conjugate to $X$, and let $L \in \Bbb Z_+$ be
such that $[-lL, lL]$ is a coding window of a topological conjugacy of $X$ onto
$\widetilde{X}$, the topological conjugacy being given by the block map
$$
\Phi: X_{[-lL, lL]} \to  \widetilde{\Sigma}.
$$

Then
$$
\widetilde {b}  = \Phi (b(db)^{4L})\in{\mathcal L}_{s-synchro}( \widetilde{X}).
$$
\end{lem}
\begin{pf}
For the proof
let $\widetilde {c} \in {\mathcal L} (\widetilde{X})$ and let $ c \in {\mathcal L}(X)  $
be such that 
$$
 \widetilde {c} = \Phi(c) .
$$
By (2.3) there is a $d^{\prime }\in \Gamma^+(c)$ such that
\begin{align*}
cd^{\prime }\in \omega^-(b). \tag {2.4}
\end{align*}
Let  $\widetilde {d}  \in
\Gamma^ {+}(  \widetilde {c}  )$ be given by
$$
\widetilde {c}\widetilde {d}\widetilde {b} = \Phi(cd^{\prime }b(db)^{4L}).
$$
(2.4) implies that 
$$
   \widetilde {c}\widetilde {d} \in \omega^{-}(  \widetilde {b} ). \qed
$$
\renewcommand{\qedsymbol}{}
\end{pf}

A topologically transitive subshift $X \subset \Sigma^{{\Bbb Z}}$ with an s-synchronizing word is called s-synchronizing (see \cite {Kr3}).
For $b \in {\mathcal L} _{s-synchro}(X)$ and $a \in \Gamma ^{+}(b)$ also $ba \in  {\mathcal L} _{s-synchro}(X)$. It follows for an s-synchronizing subshift $X$  that the set
$$
{\mathcal V} _{s-synchro}(X) = \{ \Gamma^{+}_{\infty}(b) : b \in {\mathcal L} _{s-synchro}(X) \}
$$
is the vertex set of an irreducible countable  Shannon sub-graph  of 
${\mathcal G}({\mathcal V}_{max}(X))$, that we denote by  ${\mathcal G}({\mathcal V}
_{s-synchro}(X))$, and that we call the s-synchronizing Shannon graph of $X$. ${\mathcal G}({\mathcal V}
_{s-synchro}(X))$
 presents $X$.
It follows from Lemma 4 that s-synchronization is an invariant of topological
conjugacy. The presentation of an s-synchronizing shift by its s-synchronizing
Shannon graph
is canonical. This is the content of the following theorem, 
that is shown in the same way as Theorem 1,
Lemma 4 taking the place of Lemma  3.

\begin{thm 2}
Let $X \subset \Sigma^{(\Bbb Z)}, \widetilde{X} \subset \widetilde{\Sigma}^{(\Bbb
Z)}$  be s-synchronizing subshifts, and let $\varphi: X \to \widetilde{X}  $ be a
topological conjugacy. Then
$$
\varphi ({\mathcal V}_{s-synchro}(X)) =   {\mathcal V}_{s-synchro}(\widetilde{X} ).
$$
\end{thm 2}

\subsection{a-synchronization}
For a subshift $X\subset \Sigma^{{\Bbb Z}}$, we denote for $b \in  {\mathcal L}(X)$ by $ \omega^{+}_\circ(b)$
the set of words $c \in  \omega^{+}(b)$ that appear as prefixes of sequences in $  \omega^{+}_\infty(b)$,
$$
 \omega^{+}_\circ(b) = \{x^+_{[1, n]}:  x^+\in  \omega^{+}_\infty(b), n \in{ \Bbb N}\}.
$$

A word $b$ that is admissible for a topologically transitive subshift
$X\subset \Sigma^{{\Bbb Z}}$ is called an a-synchronizing word for $X$ if $b$
satisfies the following two conditions a-s(1) and a-s(2):
\begin{cond 1}
For  $c \in{\mathcal L} (X)$ there exists a $d \in \Gamma^{-} (c)$ such 
that $d c \in
\omega^{+}(b).$
\end{cond 1}
\begin{cond 2}
For  $c \in \omega^{+}_\circ(b)$ there exists a $d \in \Gamma^{+} (c)$
such that
$cd \in\Gamma^{-} (b)$ and such that $ cdb \in \omega^{+}(b)$ and
$\omega^{+}_{\infty}(bcdb)=\omega^{+}_{\infty}(b)$.
\end{cond 2}

Condition a-s(1) is equivalent to requiring that the word $b$ is an
s-synchronizing word for the inverse of the subshift $X\subset \Sigma^{{\Bbb Z}}$.
We denote the
set of a-synchronizing words  of a topologically transitive subshift $X \subset
\Sigma^{{\Bbb Z}}$ by ${\mathcal L} _{a-synchro}(X).$ 

\begin{lem}
Let $b$ be an a-synchronizing word of the topologically transitive subshift
$X\subset \Sigma^{{\Bbb Z}}$. Then there exists a word $d \in \Gamma^{-}(b) $
such that  $db \in \omega^{+}(b)$ and $\omega^{+} _{\infty} (bdb) = \omega^{+} _{\infty} (b)$.
\end{lem}
\begin{pf}
First choose by a-s(1) a word $d^{\prime} \in  \Gamma^{-}(b) $ such that
$d^{\prime}b \in  \omega^{+}(b)$. Then choose by a-s(2) a word
$d^{\prime\prime}\in  \Gamma^{+}(d^{\prime}b)\cap   \Gamma^{-}(b)$ such that
$d^{\prime}bd^{\prime\prime}b \in \omega^{+}(b)$ and such that 
$\omega^{+}_{\infty} (bd^{\prime}bd^{\prime\prime}b) =  \omega^{+} _{\infty} (b)$. Set $ d
=  d^{\prime}bd^{\prime\prime}$.
\end{pf}
From the preceding lemma it is seen that the formulation of the following
lemma is meaningful.

\begin{lem}
Let $X 
\subset \Sigma^{(\Bbb Z)}$ be a topologically transitive subshift, and let 
\begin{align*}
b\in {\mathcal L}_{a-synchro} (X). \tag {2.5}
\end{align*}
Let $d \in {\mathcal L} (X)$ be
such that $db \in \omega^{+}(b)$ and  $\omega^{+} _{\infty} (bdb) = \omega^{+} _{\infty} (b)$, and let $l$ denote the length
of the word $db$. 

Let $\widetilde{X} \subset \widetilde{\Sigma}   ^{{\Bbb Z}}$ be a
subshift that is topologically conjugate to $X$, and let $L \in{{ \Bbb Z}}_+$ be
such that $[-lL, lL]$ is a coding window of a topological conjugacy of $X$ onto
$\widetilde{X}$, the topological conjugacy being given by the block map
$$
\Phi: X_{[-lL, lL]} \to  \widetilde{\Sigma}.
$$

Then
$$
\widetilde {b}  = \Phi ((bd)^{4L}b)\in{\mathcal L}_{a-synchro}( \widetilde{X}).
$$
\end{lem}
\begin{pf}
That $\widetilde{b}$ satisfies  a-s(1) is Lemma 4. 

For the proof that 
$\widetilde{b}$ also
satisfies a-s(2), let 
$\widetilde{c}\in \omega^{+}_{\circ}(\widetilde{b})  $, and denote the length of the word $\widetilde{c}$ by $m$. Choose an
$\widetilde{x}^{+} \in \omega^{+}_{\infty}(\widetilde{b})  $ 
such that
$
 \widetilde{x}^{+}_{[1, m]}=\widetilde{c},
$
and let $x^{+} \in \omega^{+}_{\infty}(b)$
be given by
$$
\Phi( (bd)^{2l}b    x^{+} ) =\widetilde{b}      \widetilde{x}^{+}   .
$$
Let
$$
c = x^{+} _{[1, m+lL]}.
$$
By (2.5), and according to a-s(2), there is a $d^\prime \in \Gamma ^+(c) \cap \Gamma^-(b)$ such that 
$
c d^\prime b \in \omega^+_\infty (b),
$
and
$
\omega^+_\infty (bc d^\prime b) = \omega^+_\infty (b).
$
A $\widetilde {d} \in \Gamma^+( \widetilde  {c} \cap \Gamma^-( \widetilde {b})$ such that
$$
\widetilde{c} \widetilde { d } \widetilde {b} \in \omega^+_\infty (  \widetilde {b}),
$$
is given by
$$
\Phi( (bd)^{2L} cd^{\prime} b(db)^{4L}) = \widetilde {b} \widetilde{c} \widetilde { d } \widetilde {b} .
$$

To show that also
$$
\omega^+_\infty ( \widetilde {b} \widetilde{c} \widetilde { d } \widetilde {b} ) \subset  \omega^+_\infty ( \widetilde {b} ),
$$
let
 \begin{align*}
  \widetilde {y }^+\in \omega^+_\infty ( \widetilde {b} \widetilde{c} \widetilde { d } \widetilde {b} ) . \tag {2.6}
\end{align*}
There is a 
$$
y^+\in \Gamma^+_\infty ((bd)^{2L}bcd^\prime b(db)^{3L})
$$
given by
$$
\Phi( b(db)^{3L}y^+)   =  \widetilde {b} \widetilde{c} \widetilde { d } \widetilde {b}   \widetilde {y }^+.
$$
(2.6) implies that
$$
y^+ \in\omega^+_\infty ( (bd)^{2L} bcd^\prime b(db)^{3L}  ),
$$
which implies that
$$
y^+ \in\omega^+_\infty ( b   ),
$$
which then implies that
$$
 \widetilde {y }^+ \in \omega^+_\infty ( \widetilde {b} ).\qed
$$
\renewcommand{\qedsymbol}{}
\end{pf}

For an a-synchronizing word $b$ of the topologically transitive subshift $X
\subset \Sigma^{{\Bbb Z}}$ and for a word $a \in \omega^{+}_\circ(b)$ the word $ba$ is
again a-synchronizing for $X$. It follows that the set 
$$
{ \mathcal V}_{a-synchro}(X) = \{\omega^{+}_{\infty}(b): b \in { \mathcal L} _{a-synchro}(X) \}
$$ 
is the vertex set of a Shannon sub-graph of ${\mathcal G}({\mathcal V}_{max}(X))$, that we call the a-synchronizing Shannon graph of $X$, and that we denote by ${ \mathcal G}  _{a-synchro} 
(X)$. ${ \mathcal G}  _{a-synchro}  (X)$ is the union of its irreducible components
all of which present $X$. One can see from condition a.s(2) and from Lemma 4 that
the number of irreducible component of ${ \mathcal G}  _{a-synchro}  (X)$ is an
invariant of topological conjugacy. In an attempt to maintain an anology with
synchronizaion and s-synchronization we say that a topologically transitive
subshift $X \subset \Sigma^{{\Bbb Z}}$ is a-synchronizing, if $X$ has a
a-synchronizing word and if ${ \mathcal G}  _{a-synchro}  (X)$ is irreducible (comp. \cite {Kr3}). The
presentation of an a-synchronizing subhift by its a-synchronizing Shannon graph
is canonical. This is the content of the following theorem, that is shown in the same way as Theorem 1,
Lemma 6 taking the place of Lemma  3. 

\begin{thm 3}
Let $X \subset \Sigma^{{\Bbb Z}}, \widetilde{X} \subset \widetilde{\Sigma}^{{\Bbb
Z}}$  be a-synchronizing subshifts, and let $\varphi: X \to \widetilde{X}  $ be a
topological conjugacy. Then
$$
\varphi ({ \mathcal V} _{a-synchro}(X)) =    { \mathcal V}_{a-synchro}(\widetilde{X} ).
$$
\end{thm 3}

\par\noindent Institut f\"ur Angewandte Mathematik
\par\noindent Universit\"at Heidelberg
\par\noindent Im Neuenheimer Feld 294, 69120 Heidelberg, Germany
\par\noindent e-mail: krieger@@math.uni-heidelberg.de
\end{document}